\documentclass{ifacconf}
\usepackage{graphicx}
\usepackage{mathrsfs,amssymb,amsmath,array,xcolor,cite,color}
\usepackage{mathbbold,bbm}
\usepackage{subfigure}
\usepackage{bm}
\usepackage[all]{xy}
\usepackage{chemarrow}
\usepackage{natbib}

 \newtheorem{example}{Example}

\begin{document}
\begin{frontmatter}

\title{On Stability of Two Kinds of Delayed Chemical Reaction Networks\thanksref{footnoteinfo}} 

\thanks[footnoteinfo]{This work was funded by the National Nature Science Foundation of China under Grant No. 12071428 and 62111530247, and the Zhejiang Provincial Natural Science Foundation of China under Grant No. LZ20A010002.}

\author[First]{Xiaoyu Zhang} 
\author[First]{Chuanhou Gao}
\author[Third]{Denis Dochain}

\address[First]{School of Mathematical Sciences, Zhejiang
University, Hangzhou 310027, China(e-mail: Xiaoyu\_Z@zju.edu.cn, gaochou@zju.edu.cn).}
\address[Third]{ICTEAM, UCLouvain, B\^{a}timent Euler, avenue Georges Lema\^{i}tre 4-6, 1348 Louvain-la-Neuve, Belgium (e-mail: denis.dochain@uclouvain.be)}

\begin{abstract}                
For the networks that are linear conjugate to complex balanced systems, the delayed version may include two classes of networks: one class is still linear conjugate to the delayed complex balanced network, the other is not. In this paper, we prove the existence of the first class of networks, and emphasize the local asymptotic stability relative to a certain defined invariant set. For the second class of systems, we define a special subclass and derive the local asymptotic stability for the subclass. Two examples are provided to illustrate our results. 


\end{abstract}

\begin{keyword}
Delayed chemical reaction network, mass-action kinetics, linear conjugate, complex balanced system, local asymptotic stability
\end{keyword}

\end{frontmatter}

\section{Introduction}
Chemical reaction networks (CRNs) are widespread in the normal operation of nature and the realization of various biological functions by organisms. Studying their dynamic properties can guide people to produce and live more scientifically.
%
The related research \citep{Feinberg1972Complex,Horn1972} can be directly applied to systems biology, synthetic biology, industrial chemistry, ecosystems, etc. They are even used in fields that seem not related to reactions, such as medicine \citep{2010AnIntroduction}, electricity \citep{samardzija1989nonlinear}, and machine learning \citep{anderson2020reaction}. Despite these facts, there are still many challenges in understanding the interaction in reactions. Typical examples include catalytic reactions whose intermediate processes are extremely complex, such as gene regulatory networks. In these cases, if traditional network models are adopted, too many variables will be involved that make it impossible to perform dynamical analysis.
A common solution is to use time delays to reduce the model complexity \citep{G2019}. There have been some application results, such as inducing gene switch in biological system \citep{Wang2012}, modeling transport system \citep{O2010}, etc.
Assuredly, the influence of time delay on the dynamical properties of the system is extremely complex. In biochemical systems, time delays are often used to produce oscillation in a stable system; conversely, they can be also used to stabilize a unstable system \citep{Intro}. 

The dynamical investigation on delayed CRNs has become active in the recent decades, such as modeling \citep{G2018}, stability analysis \citep{G2018,G20182} and persistence analysis \citep{H2019, Zhang2021}. Following these studies, this paper continues to focus on the stability issue of delayed CRNs. As one might know, in non-delayed case, complex balanced (CB) networks are locally asymptotically stable with the well-known pseudo-Helmholtz free energy as the Lyapunov function \citep{Horn1972}. The result is extended to larger set of networks that are linear conjugate to any CB system \citep{Johnston2011Linear}, referred to as $\ell c$CB systems. In delayed case, \cite{G20182} showed that delayed complex balanced (DCB) systems can maintain the local asymptotic stability of the corresponding CB systems for any time delay. This motivates us to consider whether $\ell c$CB systems have such a property too when time delays are introduced. We denote the systems under consideration by delayed $\ell c$CB (D$\ell c$CB) systems. 

There are two possibilities for the element in D$\ell c$CB systems. One is that the system is linear conjugate to a DCB system, the other is not. We prove the first class of systems are existing, labeled by $\ell c$DCB systems. Moreover, by decomposing the phase space of a $\ell c$DCB system into several equivalent invariant classes, we prove there must exist a unique equilibrium relative to the invariant set. Further, we prove that each equilibrium has local asymptotic stability relative to the defined invariant set. For the second class of systems, termed $\bar{\ell c}$DCB systems, by defining a class of special systems, named $\bar{\ell c}$DCB$_1$ systems, a subset of $\bar{\ell c}$DCB systems, we prove there must exist a unique equilibrium in each positive stoichiometric compatibility class. Finally, the local asymptotic stability of a $\bar{\ell c}$DCB$_1$ system is also captured.


This paper is organized as follows. Section \ref{II} gives preliminaries about CRNs and the delayed version. The local asymptotic stability of the $\ell c$DCB system is presented in Section \ref{III}. Section \ref{IV} derives the local asymptotic stability of the $\bar{\ell c}$DCB$_1$ system followed by an example as illustration. 
~\\~\\
\noindent{\textbf{Mathematical Notation:}}\\
\rule[1ex]{\columnwidth}{0.8pt}
\begin{description}
\item[\hspace{-0.5em}{$\mathbb{R}^n, \mathbb{R}^n_{\geq 0},\mathbb{R}^n_{>0}:$}] $n$-dimensional real space; $n$-dimensional non-negative real space; $n$-dimensional positive real space.
\item[\hspace{-0.5em}{$\bar{\mathscr{C}}_{+}, \mathscr{C}_{+}$}]: $\bar{\mathscr{C}}_{+}=C([-\tau,0];\mathbb{R}^{n}_{\geq 0}), \mathscr{C}_{+}=C([-\tau,0];\mathbb{R}^{n}_{>0})$ the non-negative, positive continuous function vectors defined on the interval $[-\tau,0]$, respectively.
\item[\hspace{-0.5em}{$x^{y_{\cdot i}}$}]: $x^{y_{\cdot i}}\triangleq\prod_{j=1}^{n}x_{j}^{y_{ji}}$, where $x,y_{\cdot i}\in\mathbb{R}^{n}$.
	\item[\hspace{-0.5em}{$\mathrm{Ln}(x)$}]: $\mathrm{Ln}(x)\triangleq\left(\ln{x_{1}}, \cdots, \ln{{x}_{n}} \right)^{\top}$, where $x\in\mathbb{R}^{n}_{>{0}}$.
\end{description}
\rule[1ex]{\columnwidth}{0.8pt}

\section{Preliminaries}\label{II}
In this section, some basic concepts about CRNs and the corresponding delayed version are given, respectively.

Consider a CRN containing $n$ chemical species, denoted by $X_1$, $X_2$,...,$X_n$, that take part in $r$ chemical reactions with the $i$-th reaction $R_i$ $(i=1,...,r)$ expressed as
\begin{equation} \label{eq:1}
R_i: \quad\sum^{n}_{j=1}y_{j,i}X_j\stackrel{}{\longrightarrow} \sum^{n}_{j=1}y'_{j,i}X_j,
\end{equation}
where the non-negative integers $y_{j,i}\in \mathbb{R}^n_{\geq 0}$ and $y'_{j,i}\in \mathbb{R}^n_{\geq 0}$ are the stoichiometric coefficient. Organize $y_{.i}=(y_{1,i},...,y_{n,i})^\top$ and $y'_{.i}=(y'_{1,i},...,y'_{n,i})^\top$, termed by complexes, then the stoichiometric subspace $\mathscr{S}$ of the network is defined by
 \begin{equation}
\mathscr{S}={\rm span}\{y'_{.i}-y_{.i}\;|\;{\rm for\; all}\;i=1,...,r\}.
\end{equation}
And $\mathscr{S}^\bot=\{a\in \mathbb{R}^n\;|\;a^\top y=0\;{\rm for\; all}\; y\in\mathscr{S}\}$ denote the orthogonal complement of the stoichiometric subspace $\mathscr{S}$.
If the mass-action kinetics is assigned to every reaction, then the reaction rate of $R_k$ will be evaluated by 
\begin{equation}
\delta_i(x)=\kappa_i\prod^{n}_{j=1}x^{y_{j,i}}_j\triangleq\kappa_ix^{y_{.i}},
\end{equation}
where $x_j\in \mathbb{R}_{\geq0}$ is the concentration of species $X_j(j=1,...,n)$,  $x=(x_1,...,x_n)^\top$ represents the state, and the positive real number $\kappa_i$ is the reaction rate constant. The dynamics of a mass-action system that captures the concentration evolution of each species is given by
\begin{equation} \label{eq:2}
\dot{x}(t)=\sum^{r}_{i=1}\kappa_ix(t)^{y_{.i}}(y'_{.i}-y_{.i}),\quad t\geq0.
\end{equation}
We usually use a quadruple $\mathcal{M}=(\mathcal{S,C,R},\bm{\kappa})$ to express a mass-action system, where $\mathcal{S},\mathcal{C},\mathcal{R},\bm{\kappa}$ are the set of species, complex, reactions and reaction rate constants, respectively. 

A positive vector $\bar{x}\in \mathbb{R}^n_{>0}$ is called a positive equilibrium of $\mathcal{M}$ if it satisfies $\dot{\bar{x}}=0$ in Eq. (\ref{eq:2}); and $\bar{x}$ is called a \emph{complex balanced} equilibrium if for any complex $\eta\in\mathbb{Z}^n_{\geq 0}$ in the network it satisfies
 \begin{equation}
\sum_{i:~y_{.i}=\eta}\kappa_i\bar{x}^{y_{.i}}=\sum_{i:~y'_{.i}=\eta}\kappa_i\bar{x}^{y_{.i}}.
\end{equation}

Inclusion of time delays in the reactions will not affect the properties related to the network structure, but affect dynamical properties a lot. A delayed mass-action system shares the same stoichiometric subspace and equilibrium with the corresponding mass-action system, but has different dynamics and non-negative stoichiometric compatibility class from the latter.  \cite{G20182,G2018} made extensive studies on delayed mass-action systems.

Now we introduce the delayed mass-action system. The dynamics of a mass-action system with time delays takes
\begin{equation} \label{eq:dde}
\dot{x}(t)=\sum^{r}_{i=1}\kappa_i[x(t-\tau_i)^{y_i}y'_i-x(t)^{y_i}y_i],\quad t\geq0.
\end{equation}
where $\tau_i\geq0$, $1\leq i\leq r$ are constant time delays. Clearly, each delayed system can be denoted as $D\mathcal{M}=(\mathcal{S,C,R},\bm{k},\bm{\tau})$ where $\mathcal{S,C,R}$ are the set of species, complex, reactions respectively and $\bm{k},\bm{\tau}$ are the vectors of reaction rate constants and time delays respectively. And the solution space of delayed system \eqref{eq:dde} is $\bar{\mathscr{C}}_{+}$. When $\tau_i=0$ holds for $1\leq i\leq r$, the system (\ref{eq:dde}) will reduce to (\ref{eq:2}). \cite{G2018} gives a equivalent class decomposition of phase space $\bar{\mathscr{C}}$ called the non-negative stoichiometric compatibility class. Each equivalent class is a forward invariant set of trajectory, i.e., the trajectory $x^{\theta}$ starting from $\theta$ always stays in the stoichiometric compatibility class $\mathcal{P}_{\theta}$ containing $\theta$.  The definition of $\mathcal{P}_{\theta}$ for the delayed system (\ref{eq:dde}) is given by
\begin{equation}\label{eq:scc}
\mathcal{P}_\theta=\{\psi\in \bar{\mathscr{C}}_+\;|\;c_a(\psi)=c_a(\theta)\; {\rm for\; all\;} a\in\mathscr{S}^\bot\},
\end{equation}
where the functional $c_a:\bar{\mathscr{C}}_+\rightarrow \mathbb{R}$ is defined by
 \begin{equation}\label{eq:zc0}
  \begin{split}
c_a(\psi)&=a^{\top}\biggl[\psi(0)+\sum^r_{i=1}\biggl(\kappa_i\int^0_{-\tau_i}\psi(s)^{y_i}ds\biggr)y_i\biggr].\\
   \end{split} 
\end{equation}
\cite{G20182} gives the Lyapunov functional of delayed complex balanced system with the following form
$V:\mathscr{C}_{+}\rightarrow \mathbb{R}_{\geq 0}$ is given by  
 \begin{equation}\small\label{eq:Vd}
\begin{split}
&V(\psi)=\sum_{j=1}^{n}(\psi_{j}(0)(\ln(\psi_{j}(0))-\ln(\bar{x}_{j})-1)+\bar{x}_{j})\\
&+\sum_{i=1}^{r}\kappa_{i}\int^{0}_{-\tau_{i}}\left\{(\psi (s))^{y_{\cdot i}}\left[\ln((\psi(s)^{y_{\cdot i}})-\ln(\bar{x}^{y_{\cdot i}})-1\right]+\bar{x}^{y_{\cdot i}}\right\} ds,
\end{split}
\end{equation}
 They also derive the existence, uniqueness, and local asymptotic stability of equilibriums of delayed complex balanced system relative to the stoichiometric compatibility class. 
\section{Stability of $\ell \MakeLowercase{c}$DCB systems}\label{III}
In this section, we study the local asymptotic stability of systems called $\ell c$DCB systems that are linear conjugate to delayed complex balanced systems relative to some invariant set.
\subsection{The existence of $\ell c$DCB systems}
Firstly, we give the definition of linear conjugacy which will be used throughout the paper.
\begin{defn}[\cite{Johnston2011Linear}]
The two mass-action systems $\mathcal{M}$ and $\mathcal{\tilde{M}}$ are linear conjugate if there exists a linear, bijective mapping: $\bm{\text{h}}:\mathbb{R}^{n}_{>0}\to\mathbb{R}^{n}_{>0}$ such that any trajectories $\Phi$ and $\tilde{\Phi}$ of $\mathcal{M},\tilde{\mathcal{M}}$ satisfy $\bm{\text{h}}(\Phi(x_0,t))=\tilde{\Phi}(\text{h}(x_0),t)$ for any $x_0\in\mathbb{R}^{n}_{>0}$.
\end{defn}
The above definition is the linear conjugacy for non-delayed systems, where $x_0\in \mathbb{R}^{n}_{>0}$. For systems with time delays, a generalized definition is that the dynamical equation of $D\mathcal{M},D\mathcal{\tilde{M}}$ satisfy that $\dot{x}=Q\dot{\tilde{x}}$. 
 Then we define the $\ell c$DCB system.
 \begin{defn}\label{def:LcDCB}
 A delayed mass-action chemical reaction system $D\mathcal{M}$ is called a $\ell c$DCB system if it is linear conjugate to a delayed complex balanced system $D\mathcal{\tilde{M}}$.
 \end{defn}
 Further we can derive the existence of $\ell c$DCB systems for each delayed complex balanced system.
 \begin{thm}\label{th:e}
 For any delayed complex balanced system $D\mathcal{\tilde{M}}=\{\mathcal{\tilde{S},\tilde{C},\tilde{R}},\bm{\tilde{\kappa}},\bm{\tilde{\tau}}\}$ and each positive diagonal matrix $Q$, there must exist the corresponding $\ell c$DCB system.
 \end{thm}
 \begin{pf}
The dynamical equation of the delayed complex balanced system $D\mathcal{\tilde{M}}$ can be described as $\dot{\tilde{x}}$. For each $Q$, finding a linear conjugate $\ell c$DCB network is equivalent to the realization of the following delayed differential equations
\begin{equation}
\begin{split}\label{eq:ddeLCB1}
\dot{x}& =\sum^{r}_{i=1}Q\tilde{\kappa}_{i}
     [(\tilde{x}(t-\tau_i))^{y_{\cdot i}}y'_{\cdot i}-(\tilde{x}(t))^{y_{\cdot i}}y_{\cdot i}]\\
     &=\sum^{r}_{i=1}Q\tilde{\kappa}_{i}\prod_{j=1}^{n}q_j^{-y_{ji}} [x(t-\tau_i)^{y_{\cdot i}}y'_{\cdot i}-x(t)^{y_{\cdot i}}y_{\cdot i}]
     \end{split}
\end{equation}
If $Q$ is a scalar matrix($q_1=\cdots=q_n=q$), let $\alpha_i=q\tilde{\kappa}_i\prod_{j=1}^{n}q_j^{-y_{ji}}=\tilde{\kappa}_iq^{1-\sum_{j=1}^{n}y_{ji}}$ and $\bm{\alpha}=(\alpha_1,\cdots,\alpha_r)$, then (\ref{eq:ddeLCB1}) can be realized as a $\ell c$DCB system $D\mathcal{M}=\{\mathcal{S,C,R},\bm{\kappa,\tau}\}$ where
$$\mathcal{S}=\mathcal{\tilde{S}};\mathcal{C}=\mathcal{\tilde{C}};\mathcal{R}=\mathcal{\tilde{R}};\bm{\tau}=\bm{\tilde{\tau}}, \bm{\kappa}=\bm{\alpha}.$$
If $Q$ is not a scalar matrix, i.e. there exists $i\neq j$, such that $q_i\neq q_j$. Without loss of generality, let $q_1=\text{max}\{q_j,j=1,\dots,n\}$. Otherwise, we can adjust the order of $\mathcal{\tilde{S}}$. Denote $\kappa_i=\tilde{\kappa_i}\prod_{j=1}^{n}q_j^{-y_{ji}}$, (\ref{eq:ddeLCB1}) can be written as
\begin{equation}\label{eq:LCBR}
    \begin{split}
        \dot{x}_{1}&=\sum_{i=1}^{r}\kappa_i q_1y'_{1i}x^{y_{.i}}(t-\tau_i)-\sum^{r}_{i=1}\kappa_i q_1y_{1i}x^{y_{.i}}(t),\\
        &~~~~~~~~~~~~~~~~~~~\cdots\\
        \dot{x}_{j}&=\sum_{i=1}^{r}\kappa_i q_jy'_{ji}x^{y_{.i}}(t-\tau_i)-\sum^{r}_{i=1}\kappa_iq_1y_{ji}x^{y_{.i}}(t)\\
        &+\sum^{r}_{i=1}\kappa_i(q_1-q_j)y_{ji}x^{y_{.i}}(t)\\
         &~~~~~~~~~~~~~~~~~~~\cdots\\
        \dot{x}_{n}&=\sum_{i=1}^{r}\kappa_i q_ny'_{ni}x^{y_{.i}}(t-\tau_i)-\sum^{r}_{i=1}\kappa_iq_1y_{ni}x^{y_{.i}}(t)\\
        &+\sum^{r}_{i=1}\kappa_i(q_1-q_n)y_{ni}x^{y_{.i}}(t)\\
    \end{split}
\end{equation}
Then we introduce one network realization of \eqref{eq:ddeLCB1}. The first two terms of the right side of the above equation can be realized as $r$ delayed reactions 
\begin{equation*}
    y_{.i}\xrightarrow{\kappa_iq_1,\tau_i} y''_{.i}, i=1, \cdots, r
\end{equation*}
where $y''_{ji}=q_j/q_1y'_{ji}$ and $\tau_i, \kappa_i, q_j, y_{.i},y'_{.i}$ share the same meaning with those in equation (\ref{eq:LCBR}). $\kappa_i q_1$ is the reaction rate constant of the $i$-th delayed reaction. The third term of the right side of equation (\ref{eq:LCBR}) are all positive because $q_1-q_j\geq 0$ for all $j=1,\cdots, n$. Thus it can be realized as several non-delayed reactions, for example
\begin{equation}
    y_{.i}\xrightarrow{\kappa_{i}} y_{.i}+\left(\begin{matrix}
    0\\(q_1-q_2)y_{2i}\\\cdots\\(q_1-q_n)y_{ni}
    \end{matrix}\right)=y'''_{.i}, i=1,\cdots,r
\end{equation}
So the (\ref{eq:LCBR}) can be realized as  $D\mathcal{M}=\{\mathcal{S},\mathcal{C},\mathcal{R},\bm{\kappa},\bm{\tau}\}$ where 
\begin{equation} 
\begin{split}
    &\mathcal{S}=\mathcal{\tilde{S}};~\mathcal{C}=\{y_{.i},y''_{.i},y'''_{.i},i=1,\dots,r\};\\
    &\mathcal{R}=\{y_{.i}\xrightarrow{\kappa_iq_1,\tau_i}y''_{.i}, y_{.i}\xrightarrow{\kappa_i}y'''_{.i}\}.
    \end{split}
\end{equation}
Thus $\{\mathcal{S},\mathcal{C},\mathcal{R},\bm{\kappa}, \bm{\tau}\}$ is a $\ell c$DCB network corresponding to the delayed complex balanced network $D\mathcal{\tilde{M}}$ and positive diagonal matrix $Q$.$\Box$
\end{pf}
\begin{rem}
From the form of dynamical equation of $\ell c$DCB in \eqref{eq:ddeLCB1}, we can obtain that the $i$-th delayed reaction in $\ell c$DCB network and that in the corresponding complex balanced network have the same reactant complex $y_{.i}$ and the delay $\tau_i$.
\end{rem}
\subsection{The stability of $\ell c$DCB systems}
This subsection derives the stability of $\ell c$DCB systems through the Lyapunov second method.
\begin{thm}\label{thm:DDE}
$D\mathcal{M}=(\mathcal{S, C,R},\bm{\kappa},\bm{\tau})$ is a $\ell c$DCB system, then all positive equilibria of the system $D\mathcal{M}$ are stable.
\end{thm}
\begin{pf}
The corresponding delayed complex balanced system of a $\ell c$DCB system can be expressed by $D\mathcal{\tilde{M}}=\{\mathcal{\tilde{S},\tilde{C},\tilde{R}},\bm{\tilde{\kappa},\tilde{\tau}}\}$.
Now we consider the following candidate Lyapunov-Krasovskii functional
$V_{L}:\mathscr{C}_{+}\rightarrow \mathbb{R}_{\geq 0}$ is given by  
 \begin{equation*}\label{eq:LCBVd}
\begin{split}\small
&V_{L}(\psi)=\sum_{j=1}^{n}q_j^{-1}(\psi_{j}(0)(\ln(\psi_{j}(0))-\ln(x^*_{j})-1)+x^*_{j})\\
&+\sum_{i=1}^{r}\tilde{\kappa}_i\prod_{j=1}^{n}q_j^{-y_{ji}}\int^{0}_{-\tau_{i}}\left\{(\psi (s))^{y_{\cdot i}}[\ln\left(\frac{\psi(s)^{y_{\cdot i}}}{x^{*y_{\cdot i}}}\right)-1]+x^{*y_{\cdot i}}\right\} ds,
\end{split}
\end{equation*}
where $q_j$ share the same meaning with that in Definition \ref{eq:ddeLCB1}. From the inequatity: for artibrary $c_1>0$, there exist $c_2>0$ that 
\begin{equation*}
    x[\ln{x}-\ln{c_1}-1]+c_1\geq c_2\ln{[1+(x-c_1)^2]}\geq 0
\end{equation*}
Thus $V_L\geq 0$ and $V_L(\psi)=0$ iff $\psi$ is a positive equilibrium of $\ell c$DCB system. The following part devoted to deriving that $V_L$ is also disspative.
\begin{equation*}\label{eq:LCBd}
\begin{split}
 &\dot{V}_{L}(x(t))=Q^{-1}\text{Ln}\left(\frac{x(t)}{x^{*}}\right)\dot{x}\\&+\sum_{i=1}^{r}\tilde{\kappa}_{i}\prod_{j=1}^{n}q_j^{-y_{ji}}x(t)^{y_{.i}}\left(\ln{\left(\left\{\frac{x(t)}{x^{*}}\right\}^{y_{.i}}\right)}-1\right)\\&-\sum_{i=1}^{r}\tilde{\kappa}_{i}\prod_{j=1}^{n}q_j^{-y_{ji}}x(t-\tau_{i})^{y_{.i}}\left(\ln{\left(\left\{\frac{x(t-\tau_i)}{x^{*}}\right\}^{y_{.i}}\right)}-1\right)\\
 &= \text{Ln}\left(\frac{\tilde{x}(t)}{\tilde{x}^{*}}\right) \dot{\tilde{x}}+\sum_{i=1}^{r}\tilde{\kappa}_i\tilde{x}(t)^{y_{.i}}\left(\ln{\left(\left\{\frac{\tilde{x}(t)}{\tilde{x}^{*}}\right\}^{y_{.i}}\right)}-1\right)\\
 \end{split}
 \end{equation*}
 \begin{equation}
 \begin{split}
 &-\sum_{i=1}^{r}\tilde{\kappa}_{i}\tilde{x}(t-\tau_{i})^{y_{.i}}\left(\ln{\left(\left\{\frac{\tilde{x}(t-\tau_i)}{\tilde{x}^{*}}\right\}^{y_{.i}}\right)}-1\right)\\&=\dot{V}(\tilde{x}(t))\leq 0.
\end{split}
\end{equation}
  The last equation of (\ref{eq:LCBd}) is obtained by using the dissipativeness of the Lyapunov-Krasovskii functional $V$ along each trajectory $\tilde{x}(t)$ of the complex balanced system $D\mathcal{\tilde{M}}$ and $\dot{V}(\tilde{x}^*)=0, \tilde{x}^*\in\mathscr{C}_{+}$ if and only if $\tilde{x}^*$ is an equilibrium of the delayed complex balanced system. Thus (\ref{eq:LCBd}) reveals that $V_{L}$ is dissipative along each trajectory $x(t)$ of the $\ell c$DCB system $D\mathcal{M}$ and
\begin{equation}
\begin{split}
\dot{V}_{L}(x^*)=0&\Longleftrightarrow \dot{V}(\tilde{x}^*)=0
    \end{split}
\end{equation}
$\tilde{x}^*$ must be a positive equilibrium of $D\mathcal{\tilde{M}}$. Thus $x^*$ is a positive equilibrium of the $\ell c$DCB system $D\mathcal{M}$. Hence, any positive equilibrium of the $\ell c$DCB system is stable.
\end{pf}
But usually the local asymptotic stability of the $\ell c$DCB system does not hold relative to the invariant class---chemical stoichiometric compatibility class $\mathcal{P}_{\theta}$ defined in (\ref{eq:scc}) because of the degenerate equilibrium points. So we re-decompose the solution space $\bar{\mathscr{C}}_{+}$ of $\ell c$DCB networks:
\begin{lem}
$D\mathcal{M}=(\mathcal{S,C,R},\bm{\kappa},\bm{\tau})$ is a $\ell c$DCB system, and its corresponding delayed complex balanced system is $D\mathcal{\tilde{M}}=(\mathcal{\tilde{S},\tilde{C},\tilde{R}},\bm{\tilde{\kappa}},\bm{\tau})$. $\mathscr{S}, \tilde{\mathscr{S}}$ are the stoichiometric subspace of $D\mathcal{M},D\mathcal{\tilde{M}}$ respectively. Then the following set is the invariant set of each trajectory of $\ell c$DCB networks.
\begin{equation}\label{eq:H}
\mathcal{H}_\theta=\{\psi\in \bar{\mathscr{C}}_+\;|\;h_a(\psi)=h_a(\theta)\; {\rm for\; all\;} a\in (Q^{-1})^{\top}\tilde{\mathscr{S}}^\bot\},
\end{equation}
where the functional $h_a:\bar{\mathscr{C}}_+\rightarrow \mathbb{R}$ is defined by
\begin{equation}\label{eq:H0}
h_a(\psi)=a^{\top}\biggl[\psi(0)+\sum^{\tilde{r}}_{i=1}\biggl(\tilde{\kappa}_i\prod_{j=1}^{n}q_j^{-y_{ji}}\int^0_{-\tau_i}\psi(s)^{y_i}ds\biggr)Qy_i\biggr]
\end{equation}
\end{lem}
\begin{pf}
Consider any trajectory $x(t)$ of a $\ell c$DCB system, then we study the change of the value of $h_{a}$ along the trajectory $x(t)$.
\begin{equation}
    \frac{\text{d}h_{a}(x(t))}{\text{d}t}=a^{\top}\biggl(\sum_{i=1}^{\tilde{r}}\tilde{\kappa}_i\tilde{x}^{y_{.i}}(t-\tau_i)Q(\tilde{y}'_{.i}-\tilde{y}_{.i})\biggr)=0
\end{equation}
Thus we conclude the result.$\Box$
\end{pf}
Now we consider the situation of equilibriums of $\ell c$DCB system relative to the above invariant set.
\begin{lem}\label{lem:1}
For a $\ell c$DCB network denoted as $D\mathcal{M}=(\mathcal{S,C,R},\bm{\kappa},\bm{\tau})$, the positive equilibrium in each invariant set $\mathcal{H}_{\theta}$ of $D\mathcal{M}$ defined as (\ref{eq:H}) is unique.
\end{lem}
\begin{pf}
$\mathcal{H}_{\theta}$ defined in \eqref{eq:H} is a arbitrary invariant set of $D\mathcal{M}$. The corresponding delayed complex balanced system of $D\mathcal{M}$ is denoted as $D\mathcal{\tilde{M}}=(\mathcal{\tilde{S},\tilde{C},\tilde{R}},\bm{\tilde{\kappa},\bm{\tilde{\tau}}})$ and $\mathcal{P}_{\tilde{\theta}}$ is the chemical stoichiometric compatibility class of $D\mathcal{\tilde{M}}$ containing $\tilde{\theta}$ where $\theta=Q\tilde{\theta}$. Denoting the unique positive equilibrium of $\mathcal{P}_{\tilde{\theta}}$ as $\tilde{x}^*$, we claim that $x^*=Q\tilde{x}^*$ must be the unique positive equilibrium of the invariant set $\mathcal{H}_{\theta}$. In order to derive this result, we just need to verify the values of $h_{a}(\theta)$, $h_{a}(x^*)$ for each $a$. $\tilde{x}^*$ is in $\mathcal{P}_{\tilde{\theta}}$, thus $c_{\tilde{a}}(\tilde{x}^*)=c_{\tilde{a}}(\tilde{\theta})$ where $\tilde{a}\in\tilde{\mathscr{S}}^{\bot}$. Also,
\begin{equation*}
\begin{split}
   h_{a}(\theta)&=\tilde{a}^{\top}Q^{-1}\biggl[Q\theta(0)+\sum_{i=1}^{r}\left(\tilde{\kappa}_i\int^{0}_{-\tau_i}\tilde{\theta}(s)^{y_{.i}}ds\right)Qy_{.i}\biggr]\\
   &=c_{\tilde{a}}(\tilde{\theta})
   \end{split}
\end{equation*}
and $h_{a}(x^*)=c_{\tilde{a}}(\tilde{x}^*)$. Thus $h_{a}(x^*)=h_{a}(\theta)$, i.e., $x^*$ is in the invariant set $h_{a}(\theta)$. Further, $x^*=D\tilde{x}^*$ is the positive equilibrium of $h_{a}(\theta)$. Similarly, if $\mathcal{H}_(\theta)$ has another positive equilibrium except $x^*$, the stoichiometric compatibility class $\mathcal{P}_{\tilde{\theta}}$ of system $D\mathcal{\tilde{M}}$ also has more than one equilibrium. This is obviously contradiction to the uniqueness of positive equilibrium of the delayed complex balanced system. Thus we derive the uniqueness and existence of equilibrium in $\ell c$DCB systems relative to $\mathcal{H}_{\theta}$.$\Box$
\end{pf}
\begin{thm}\label{th:lasldcb}
$D\mathcal{M}=(\mathcal{S,C,R},\bm{\kappa},\bm{\tau})$ is a $\ell c$DCB network, then each positive equilibrium $x^*$ is local asymptotic stability relative to the invariant set $\mathcal{H}_{\theta}$ containing $x^*$.
\end{thm}
\begin{pf}
This result can be derived from Theorem \ref{thm:DDE}, Lemma \ref{lem:1}.$\Box$
\end{pf}
\begin{example}
Consider the following delayed complex balanced system $D\mathcal{\tilde{M}}$:
	\begin{align*}
		\xymatrix{A\ar ^{\tilde{\tau}_1,~\tilde{k}_1~~~} [r] & 2B \ar ^{~\tilde{\tau}_2,~\tilde{k}_2} [d]\\
			& 2A+2B\ar ^{~\tilde{\tau}_3,~\tilde{k}_3~~} [lu]}
	\end{align*}
	By choosing $\tilde{k}_1=\tilde{k}_2=\tilde{k}_3=1$, the dynamical equation of $D\mathcal{\tilde{M}}$ can be written as:
	\begin{equation}
	\begin{split}
	\dot{\tilde{x}}_A&=-\tilde{x}_A(t)+2\tilde{x}_B^2(t-\tilde{\tau}_2)+\tilde{x}_A^2\tilde{x}_B^2(t-\tilde{\tau}_3)-2\tilde{x}_A^2\tilde{x}_B^2(t);\\
	\dot{\tilde{x}}_B&=2\tilde{x}_A(t-\tilde{\tau}_{1})+2\tilde{x}_B^2(t-\tilde{\tau}_2)-2\tilde{x}^2_{B}(t)-2\tilde{x}_A^2\tilde{x}_B^2(t)
	\end{split}
\end{equation}
If $Q=\text{diag}(2,2)$ is a scalar matrix, the dynamical equation of the corresponding $\ell c$DCB should be:
	\begin{equation}
		\begin{split}
			\dot{x}_A&=-x_A(t)+x_B^2(t-\tilde{\tau}_2)+\frac{1}{8}x_A^2x_B^2(t-\tilde{\tau}_3)-\frac{1}{4}x_A^2x_B^2(t)\\
			\dot{x}_B&=2x_A(t-\tilde{\tau}_{1})+x_B^2(t-\tilde{\tau}_2)-x^2_{B}(t)-\frac{1}{4}\tilde{x}_A^2\tilde{x}_B^2(t)\\
		\end{split}
	\end{equation}
	Then one $\ell c$DCB realization as described in Theorem \ref{th:e} is 
		\begin{align*}
	\xymatrix{A\ar ^{\tau_1,~k_1~~~} [r] & 2B \ar ^{~\tau_2,~k_2} [d]\\
		& 2A+2B\ar ^{~\tau_3,~k_3~~} [lu]}
\end{align*}
where $k_1=1,k_2=\frac{1}{2},k_3=\frac{1}{8}$ and $\tau_i=\tilde{\tau}_i$.

If $Q=\text{diag}(2,1)$ is not a scalar matrix, the dynamical equation of the $\ell c$DCB determined by $D\mathcal{\tilde{M}}$ and $Q$ is 
	\begin{equation}\label{ex:L-DDE}
	\begin{split}
		\dot{x}_A&=-x_A(t)+4x_B^2(t-\tilde{\tau}_2)+\frac{1}{2}x_A^2x_B^2(t-\tilde{\tau}_3)-x_A^2x_B^2(t)\\
		\dot{x}_B&=x_A(t-\tilde{\tau}_{1})+2x_B^2(t-\tilde{\tau}_2)-2x^2_{B}(t)-\frac{1}{2}\tilde{x}_A^2\tilde{x}_B^2(t)\\
	\end{split}
\end{equation}
The following delayed system is the corresponding $\ell c$DCB realization shown in Theorem \ref{th:e}
	\begin{equation}
	\begin{split}
	A&\xrightarrow{k_1,\tau_1}B\\
	2B&\xrightarrow{k_2,\tau_2}2A+B\\
	2A+2B&\xrightarrow{k_3,\tau_3}A\\
	2B&\xrightarrow{k_4}4B\\
	2A+2B&\xrightarrow{k_5}2A+4B
	\end{split}
\end{equation}
where $k_1=1,k_2=2,k_3=\frac{1}{2},k_4=1,k_5=\frac{1}{4}$, $\tau_i=\tilde{\tau}_i, i=1,2,3$ and $\tau_4=0,\tau_5=0$.
Note that the above realization is not unique, the dynamical equation of the following delayed $\ell c$DCB system is also \eqref{ex:L-DDE}.
	\begin{equation}
	\begin{split}
		A&\xrightarrow{k_1,\tau_1}B\\
		2B&\xrightarrow{k_2,\tau_2}4A+2B\\
		2A+2B&\xrightarrow{k_3,\tau_3}2A\\
		2A+2B&\xrightarrow{k_4}A+2B
	\end{split}
\end{equation}
where $k_1=1, k_2=1,k_3=\frac{1}{4},k_4=\frac{1}{2}$,  $\tau_i=\tilde{\tau}_i, i=1,2,3$ and $\tau_4=0$.

Thus from Theorem \ref{th:lasldcb}, above delayed $\ell c$DCB systems are all local asymptotic stability relative to the invariant class.
\end{example}

\section{Stability of $\bar{\ell \MakeLowercase{c}}$DCB$_1$ systems}\label{IV}
 When time delays are introduced to non-delayed systems $\mathcal{M}$ linear conjugate to complex balanced systems $\tilde{\mathcal{M}}$, the delayed version of some system $\mathcal{M}$ denoted as $D\mathcal{M}$ is no longer conjugate to the delayed complex balanced system $D\mathcal{\tilde{M}}$. In this section, we study the local asymptotic stability of a special case of $D\mathcal{M}$ called the $\bar{\ell c}$DCB$_1$ system.
\subsection{Problem statement}\label{pm}
We illustrate our motivation through the following two delayed systems
\begin{equation}\label{eq:N1}
\begin{split}
    &D\mathcal{M}: 3A\xrightarrow{1, \tau_1} A+2B\xrightarrow{2,\tau_2} 2A+B\\
 &D\mathcal{\tilde{M}}: \xymatrix{
	3A\ar @{ -^{>}}^{1,\tau_1}  @< 1pt> [r]& A+2B \ar  @{ -^{>}}^{1,\tau_2}  @< 1pt> [l] }.
\end{split}
	\end{equation}
The above two networks with no time delay ($\tau_i=0, i=1,2$) denoted as $\mathcal{M},\mathcal{\tilde{M}}$ share the same dynamical equation:
\begin{align}
\left\{
	\begin{array}{ll}
		\dot{x}_A(t)=-2x_A^3+2x_Ax_B^2,
		\\
		\dot{x}_B(t)=2x_A^3-2x_Ax_B^2.
	\end{array}
	\right.	
\end{align}

 Note that the network structure of system $\mathcal{\tilde{M}}$ is weakly reversible and zero deficiency, thus $\mathcal{\tilde{M}}$ is complex balanced system \citep{Feinberg1972Complex}. The local asymptotic stability of $\mathcal{\tilde{M}}$ can be derived by the Lyapunov function--- Pseudo-Helmholtz free energy function. Although $\mathcal{M}$ is not a complex balanced system even not a weakly reversible system, but the local asymptotic stability of $\mathcal{M}$ can be also derived using the dynamic equivalence.
 
 But when time delays are introduced into reactions, the situation will be completely different. The delayed systems $D\mathcal{M}$ and $D\mathcal{\tilde{M}}$ have completely different dynamics and are not linear conjugate to each other.
The dynamical equation of delayed system $D\mathcal{M}$ can be written as
\begin{align}\label{eq:N1de}
\left\{
	\begin{array}{ll}\dot{x}_A(t)&=x_A^3(t-\tau_1)-3x_A^3(t)+4x_Ax_B^2(t-\tau_2)-2x_Ax_B^2(t)\\
	    \dot{x}_B(t)&=2x_A^3(t-\tau_1)+2x_Ax_B^2(t-\tau_2)-4x_Ax_B^2(t);
	\end{array}
	\right.	
\end{align}
But the dynamics of delayed complex balanced network $D\mathcal{\tilde{M}}$ is
\begin{align*}
\left\{
	\begin{array}{ll}
	\dot{x}_A(t)&=x_A^3(t-\tau_1)-3x_A^3(t)+3x_Ax_B^2(t-\tau_2)-x_Ax_B^2(t)\\
	    \dot{x}_B(t)&=2x_A^3(t-\tau_1)-2x_Ax_B^2(t);
	\end{array}
	\right.	
\end{align*}
$D\mathcal{M}$ is not a $\ell c$DCB network that we introduced in section \ref{III}. In this section, we focus on the local asymptotic stability of this kind of delayed networks like $D\mathcal{M}$ which are not $\ell c$DCB networks but its non-delayed system $\mathcal{M}$ shares the same dynamical equation with a complex balanced network $\mathcal{\tilde{M}}$.
\subsection{Local asympototic stability of $\bar{\ell c}$DCB$_1$ network}
This subsection obtains the local asymptotic stability of one type of delayed systems which is a generalization of the system $D\mathcal{M}$ proposed in subsection \ref{pm}.
\begin{defn}\label{def:N1}
A delayed system $D\mathcal{M}=(\mathcal{S},\mathcal{C},\mathcal{R}, \bm{\kappa}, \bm{\tau})$ is called a $\bar{\ell c}$DCB$_1$ system if there exists a complex balanced mass-action system $D\mathcal{\tilde{M}}=(\mathcal{\tilde{S}},\mathcal{\tilde{C}},\mathcal{\tilde{R}},\bm{\tilde{\kappa}},\bm{\tilde{\tau}})$ such that
\begin{itemize} 
\item $\mathcal{S}=\mathcal{\tilde{S}}$; 
\item $\text{card}(\mathcal{\tilde{R}})=\text{card}(\mathcal{R})$; And for $i$-th reaction of two systems, there exist $y_{.i}=\tilde{y}_{.i}$ and the reaction vectors satisfiy $v_{.i}=b_{i}\tilde{v}_{.i}$, where $b_{i}=\frac{\tilde{\kappa}_{i}}{\kappa_{i}}\leq 1$ is a positive constant.
\end{itemize}
\end{defn}
\begin{rem}
Consider a $\bar{\ell c}$DCB$_1$ system $D\mathcal{M}=\{\mathcal{S,C,R},\bm{\kappa,\tau}\}$ with
$$\mathcal{R}=\{y_{.i}\xrightarrow{\kappa_i,\tau_i}y'_{.i}\}=\{y_{.i}\xrightarrow{\kappa_{i},\tau_i}y_{.i}+v_{i}\}.$$
The corresponding complex balanced network of $D\mathcal{M}$ is $D\mathcal{\tilde{M}}=(\mathcal{\tilde{S}},\mathcal{\tilde{C}},\mathcal{\tilde{R}},\tilde{\bm{\kappa}},\tilde{\bm{\tau}})$. $\mathcal{\tilde{R}}$ in system $D\mathcal{\tilde{M}}$ is expressed by: $$\mathcal{\tilde{R}}=\{\tilde{y}_{.i}\xrightarrow{\tilde{\kappa}_{i},\tilde{\tau}_i}\tilde{y}'_{.i}\}=\{y_{.i}\xrightarrow{\tilde{\kappa}_{i},\tilde{\tau}_i}y_{.i}+\tilde{v}_{.i}\}.$$

 The dynamical equation of the delayed system $D\mathcal{M}$ 
\begin{equation*}
    D\mathcal{M}: \dot{x}(t)=\sum^{r}_{i=1}\kappa_i
     [(x(t-\tau_i))^{y_{\cdot i}}(y_{.i}+b_{i}\tilde{v}_{i})-(x(t))^{y_{\cdot i}}y_{\cdot i}]
     \end{equation*}
     is different from the dynamical equation of $D\mathcal{\tilde{M}}$ even if choosing  $\tilde{\tau}_{i}=\tau_i$ for each $i=1,\cdots,r$. 
     \begin{equation*}
          D\mathcal{\tilde{M}}:\dot{x}(t)=\sum^{r}_{i=1}\tilde{\kappa}_i
     [(x(t-\tau_i))^{y_{\cdot i}}(y_{.i}+\tilde{v}_{i})-(x(t))^{y_{\cdot i}}y_{\cdot i}].
     \end{equation*}
     Also, two networks are not linear conjugate.
\end{rem}
Now we give the local asymptotic stability of the delayed system $\bar{\ell c}$DCB$_1$ by using the Lyapunov second method. 
\begin{lem}[existence and uniqueness]
For a $\bar{\ell c}$DCB$_1$ system $D\mathcal{M}=(\mathcal{S},\mathcal{C},\mathcal{R},\bm{\kappa}, \bm{\tau})$, each stoichiometric compatibility class of $D\mathcal{M}$ contains a unique positive equilibrium.
\end{lem}
\begin{pf}
$\mathcal{P}_{\theta}$ defined as (\ref{eq:scc}) denotes arbitrary stoichiometric compatibility class of the $\bar{\ell c}$DCB$_1$ system $D\mathcal{M}$. $\bm{a}=\{a_1,\cdots,a_{n-s}\}$ is a set of basis of $\mathscr{S}^{\bot}$ where $s$ denotes the dimension of $\mathscr{S}$. If $n=s$, this lemma holds obviously. If $s<n$, let $c_{a_i}(\theta)=M_i$. The stoichiometric compatibility class $\mathcal{P}_{\theta}$  can be expressed as 
\begin{equation*}
    \mathcal{P}_\theta=\{\psi\in \bar{\mathscr{C}}_+\;|\;c_{a_{i}}(\psi)=M_i\; {\rm for\; all\;} a_i\in\bm{a}\},
\end{equation*}
Now we consider the delayed complex balanced system $D\mathcal{\tilde{M}}=(\mathcal{\tilde{S},\tilde{C},\tilde{R}},\bm{\tilde{\kappa}},\bm{\tilde{\tau}})$, 
where the time delay of $i$-th reaction $\tilde{\tau}_i=\tau_i /b_i$ and $\mathcal{\tilde{S},\tilde{C},\tilde{R}},\bm{\kappa}$ share the same meaning with Definition \ref{def:N1}. Also, Definition \ref{def:N1} reveals that stoichiometric subspaces of systems
$D\mathcal{M}$ and $D\mathcal{\tilde{M}}$ are the same, thus $\mathscr{S}^{\bot}$ are also same. In this case, for the same conserved quantity $M_i,i=1,\cdots,n-s$, the corresponding stoichiometric compatibility class $\mathcal{P}_{\theta}$ of the system $D\mathcal{M}$ and $\mathcal{\tilde{P}}_{\theta}$ of the system $D\mathcal{\tilde{M}}$ contain same constant points $\psi\in \mathscr{C}_{+}$, which can be seen that
\begin{equation*}
  c_{a_i}(\psi)=\psi+\sum_{i=1}^{r}\kappa_i\psi^{y_{.i}}\tau_i y_{.i}=\psi+\sum_{i=1}^{r}\tilde{\kappa}_i/b_i\psi^{y_{.i}}\tau_i y_{.i}\\
  \end{equation*}
  \begin{equation*}
  =\psi+\sum_{i=1}^{r}\tilde{k}_i\psi^{y_{.i}}\tilde{\tau_i} y_{.i}=\tilde{c}_{a_i}(\psi)=M_i
\end{equation*}
From the fact that each positive stoichiometric compatibility class in the delayed complex balanced system $D\mathcal{\tilde{M}}$ contains and only contains one equilibrium. We can conclude the existence and uniqueness of positive equilibrium in a $\bar{\ell c}$DCB$_1$ network.$\Box$
\end{pf}
\begin{figure*}
\centering
\subfigure[$\tau=(0.1,1), x_0(s)=(5,1) $]{\includegraphics[height=5.5cm,width=8.8cm]{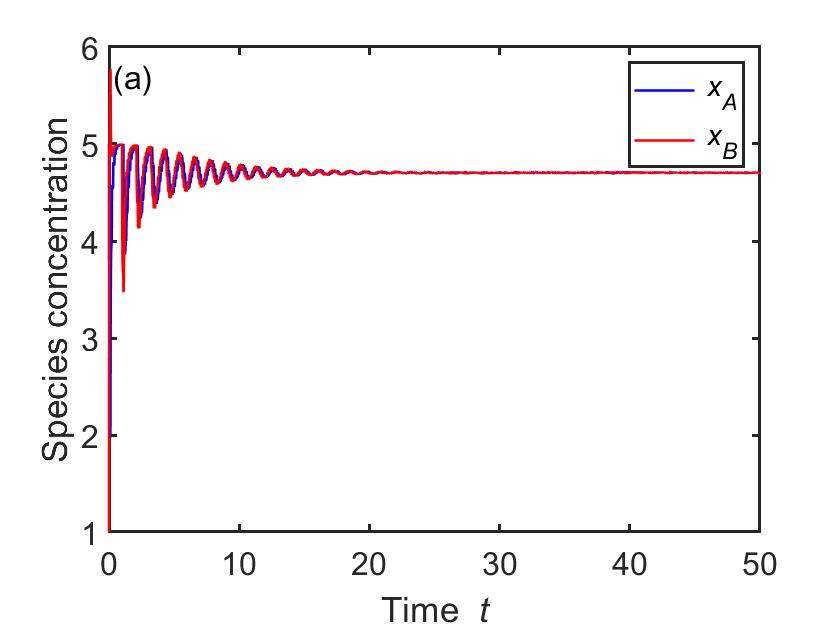}}
\subfigure[$\tau=(2,0.5), x_0(s)=(5,1)$]{\includegraphics[height=5.5cm,width=8.8cm]{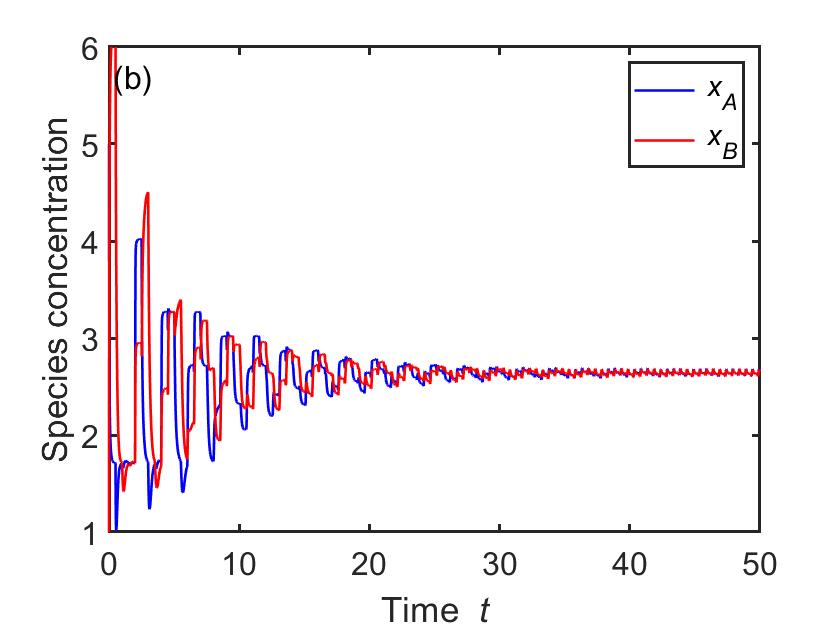}}
\quad
\subfigure[$\tau=(0.1,1), x_0(s)=(\sin{s}+1,\cos{s}+1)$]{\includegraphics[height=5.5cm,width=8.8cm]{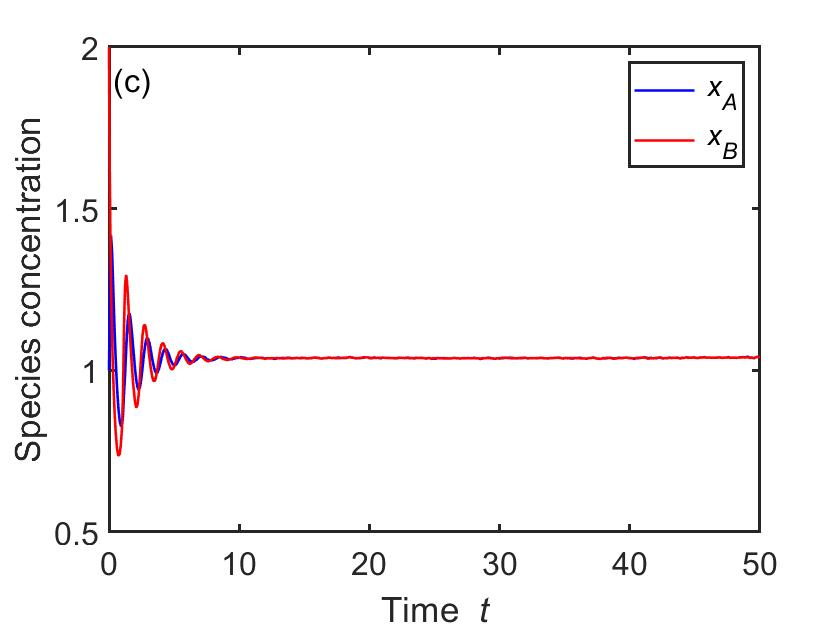}}
\subfigure[$\tau=(2,0.5), x_0(s)=(\sin{s}+1,\cos{s}+1)$]{\includegraphics[height=5.5cm,width=8.8cm]{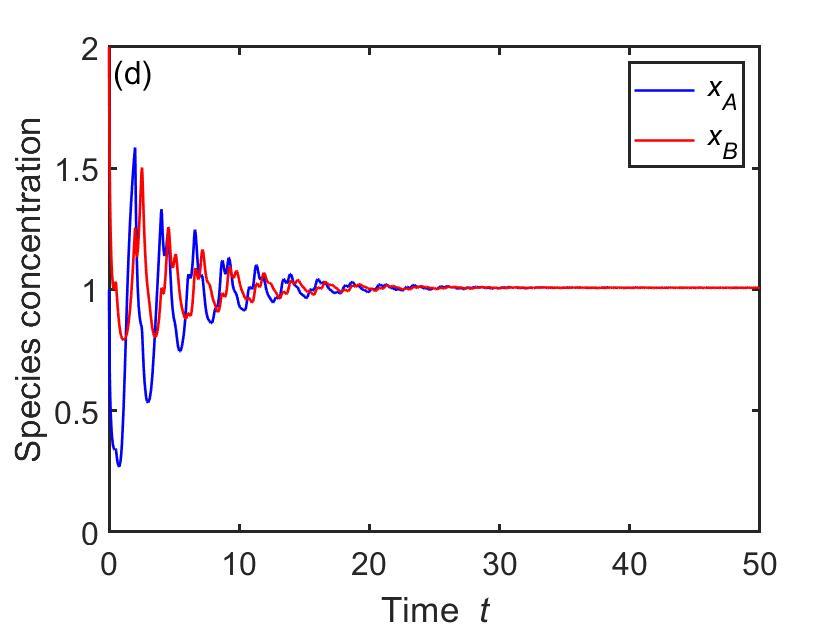}}
 \caption{The evolution behavior of Example \ref{ex:1} with different initial datas and different time delays.}\label{fig:1}
\end{figure*}
\begin{thm}\label{thm:slcb}
If $D\mathcal{M}=(\mathcal{S},\mathcal{C},\mathcal{R},\bm{\kappa}, \bm{\tau})$ is a delayed $\bar{\ell c}$DCB$_1$ network, each positive equilibrium of $D\mathcal{M}$ is local asymptotic stability.
\end{thm}
\begin{pf}
The following part is devoted to proving that the pesudo-Helmholtz functional proposed in (\ref{eq:Vd}) also can be the Lyapunov functional of the system $D\mathcal{M}$ denoted by $V_D(x(t))$.
The $\dot{V}_{D}(x(t)$ can be divided into two parts($\dot{V}_{D}(x(t))=A+B$):  
 \begin{equation*}
   \begin{split}
       A=&\text{Ln}(\frac{x(t)}{\bar{x}})\sum_{i=1}^{r}\tilde{\kappa}_i[x(t-\tau_i)^{y_{.i}}(y_{.i}+\tilde{v}_{.i})-x(t)^{y_{.i}}y_{.i}]\\+&\sum_{i=1}^{r}\tilde{\kappa}_ix(t)^{y_{.i}}[\ln{\frac{x(t)}{\bar{x}}}^{y_{i}}-1]\\
       &-\sum^{r}_{i=1}\tilde{\kappa}_{i}x(t-\tau_i)^{y_{.i}}[\ln{\left(\frac{x(t-\tau_i)^{y_{.i}}}{\bar{x}}\right)}-1]\\
      \end{split}
      \end{equation*}
       \begin{equation*}
   \begin{split}
       B=&\text{Ln}(\frac{x(t)}{\bar{x}})\sum_{i=1}^{r}(\kappa_i-\tilde{\kappa}_{i})[x(t-\tau_i)^{y_{.i}}y_{.i}-x(t)^{y_{.i}}y_{.i}]\\
       +&\sum_{i=1}^{r}(\kappa_i-\tilde{\kappa}_{i})x(t)^{y_{.i}}[\ln{\left(\frac{x(t)}{\bar{x}}\right)^{y_{.i}}}-1]
       \\&-\sum^{r}_{i=1}(\kappa_{i}-\tilde{\kappa}_{i})x(t-\tau_i)^{y_{.i}}[\ln{\left(\frac{x(t-\tau_i)}{\bar{x}}\right)^{y_{.i}}}-1]
   \end{split}
\end{equation*}
Actually, $A$ is the derivative of the Lyapunov functional of the delayed complex balanced system $D\mathcal{\tilde{M}}$ with respect to time. Therefore, from the dissipation of the delayed complex balanced network, we know that $A\leq 0$ and equality holds iff $x$ is a equilibrium of the $D\mathcal{M}$ and $D\mathcal{\tilde{M}}$. 

Now we consider the part B, let $\bar{\kappa}_{i}=\kappa_i-\tilde{\kappa}_i$ 
\begin{equation*}
\begin{split}
    B&=\sum^{r}_{i=1}\bar{\kappa}_{i}\left[x(t-\tau_i)^{y_{.i}}\left(\ln{\left(\frac{x(t)}{\bar{x}}\right)^{y_{.i}}}-\ln{\left(\frac{x(t-\tau_i)}{\bar{x}}\right)^{y_{.i}}}\right)\right]\\
    &-\sum^{r}_{i=1}\bar{\kappa}_{i}\left[x(t)^{y_{.i}}-x(t-\tau_i)^{y_{.i}}\right]
    \end{split}
 \end{equation*}
From the inequation $e^a(b-a)\leq e^b-e^a$ (equality holds iff $a=b$), for each $i$
\begin{equation*}
\begin{split}
    &x(t-\tau_i)^{y_{.i}}\left(\ln{\left(\frac{x(t)}{\bar{x}}\right)^{y_{.i}}}-\ln{\left(\frac{x(t-\tau_i)}{\bar{x}}\right)^{y_{.i}}}\right)\\
    &\leq x(t)^{y_{.i}}-x(t-\tau_i)^{y_{.i}}.
    \end{split}
\end{equation*}
Thus combining above two equations, it is obviously that the $B\leq 0$ with equality iff $x(t)=x(t-\tau_i)$ for any $\tau_i$. 

Consequently, $\dot{V}_D(x(t))=A+B\leq 0$. And $\dot{V}_D(x^*)=0$ iff $x^*$ is a constant positive equilibrium. Thus each equilibrium of the system $D\mathcal{M}$ is locally asymptotically stable.$\Box$
\end{pf}

\begin{example}\label{ex:1}
Now we reconsider the system $D\mathcal{M}$ in the equation (\ref{eq:N1}), $D\mathcal{\tilde{M}}$ is the corresponding complex balanced network by taking $$b_1=1, b_2=\frac{1}{2}; v_{.1}=\tilde{v}_{.1}=(-2,2)^{T}, v_{.2}=\frac{1}{2}\tilde{v}_{.2}=(1,-1)^{T}.$$
 Thus from theorem \ref{thm:slcb}, each positive equilibrium of $D\mathcal{M}$ is locally asymptotically stable. Through the numerical analysis of the $D\mathcal{M}$ with different time delays and initial points, the local asymptotic stability of the equilibrium point is further intuitively obtained(See Fig. \ref{fig:1}).
\end{example}
\bibliography{ifacconf}             
                                                   







\end{document}